# ON HOEFFDING'S INEQUALITIES[1]

### By Vidmantas Bentkus

*Vilnius Institute of Mathematics and Informatics,*
*and Vilnius Pedagogical University*

In a celebrated work by Hoeffding [*J. Amer. Statist. Assoc.* **58** (1963) 13–30], several inequalities for tail probabilities of sums $M_n = X_1 + \cdots + X_n$ of bounded independent random variables $X_j$ were proved. These inequalities had a considerable impact on the development of probability and statistics, and remained unimproved until 1995 when Talagrand [*Inst. Hautes Études Sci. Publ. Math.* **81** (1995a) 73–205] inserted certain missing factors in the bounds of two theorems. By similar factors, a third theorem was refined by Pinelis [*Progress in Probability* **43** (1998) 257–314] and refined (and extended) by me. In this article, I introduce a new type of inequality. Namely, I show that $\mathbb{P}\{M_n \geq x\} \leq c\mathbb{P}\{S_n \geq x\}$, where $c$ is an absolute constant and $S_n = \varepsilon_1 + \cdots + \varepsilon_n$ is a sum of independent identically distributed Bernoulli random variables (a random variable is called Bernoulli if it assumes at most two values). The inequality holds for those $x \in \mathbb{R}$ where the survival function $x \mapsto \mathbb{P}\{S_n \geq x\}$ has a jump down. For the remaining $x$ the inequality still holds provided that the function between the adjacent jump points is interpolated linearly or log-linearly. If it is necessary, to estimate $\mathbb{P}\{S_n \geq x\}$ *special* bounds can be used for binomial probabilities. The results extend to martingales with bounded differences. It is apparent that Theorem 1.1 of this article is the most important. The inequalities have applications to measure concentration, leading to results of the type where, up to an absolute constant, the measure concentration is dominated by the concentration in a simplest appropriate model, such results will be considered elsewhere.

**1. Introduction and results.** To illustrate the flavor of the inequalities provided below, let us start with the special case of a sum $Z_n = Y_1 + \cdots + Y_n$

Received January 2002; revised May 2003.
[1]Supported by Max Planck Institute for Mathematics, Bonn.

of bounded independent random variables such that $\mathbb{P}\{0 \leq Y_k \leq 1\} = 1$ and $\mathbf{E}X_k = p_k$ for all $k$. Then

$$(1.1) \qquad \mathbb{P}\{Z_n \geq x\} \leq e\mathbb{P}\{\varepsilon_1 + \cdots + \varepsilon_n \geq x\}, \qquad e = 2.718\ldots,$$

for integer $x \in \mathbb{Z}$, where $\varepsilon_1, \ldots, \varepsilon_n$ are independent identically distributed (i.i.d.) Bernoulli random variables that assume values 0 and 1 such that $\mathbb{P}\{\varepsilon_k = 1\} = p$ with $p = (p_1 + \cdots + p_n)/n$. The bound (1.1) is a very special case of Theorem 1.2. The following bound (1.2) is independent of $n$ and is much rougher than (1.1). Furthermore, usually bounds of type (1.1) are more convenient in applications than bounds of type (1.2). We have

$$(1.2) \qquad \mathbb{P}\{Z_n \geq x\} \leq \frac{e^3}{2}\mathbb{P}\left\{\eta \geq \frac{x}{1-p}\right\}$$

for $x$ such that $x/(1-p)$ is an integer, where $\eta$ is a Poisson random variable with parameter $\lambda$ such that

$$\lambda = pn/(1-p), \qquad \mathbb{P}\{\eta = k\} = \lambda^k \exp\{-\lambda\}/k! \qquad \text{for } k = 0, 1, 2, \ldots.$$

The Introduction is organized as follows. First formulations of the results, namely, of Theorems 1.1–1.3 are provided. Then their relationships to Hoeffding's inequalities are discussed, references are provided and the methods are explained. Theorem 1.1 seems to be the most important. It has nice applications to the measure concentration; such applications will be addressed elsewhere.

Henceforth replace the independence assumption by a martingale type dependence. Let

$$\mathcal{F}_0 = \varnothing \subset \mathcal{F}_1 \subset \cdots \subset \mathcal{F}_n \subset \mathcal{F}$$

be a family of $\sigma$ algebras of a measurable space $(\Omega, \mathcal{F})$. Let $M_n = X_1 + \cdots + X_n$ be a martingale with differences $X_k = M_k - M_{k-1}$. Define $M_0 = 0$.

The simplest thinkable nontrivial martingale is a sum $S_n = \varepsilon_1 + \cdots + \varepsilon_n$ of $n$ i.i.d. Bernoulli random variables. A random variable (or its distribution) is called Bernoulli if it assumes at most two values with positive probability. Let $\sigma > 0$ and $b > 0$. By $\varepsilon = \varepsilon(\sigma^2, b)$ denote a Bernoulli random variable such that

$$(1.3) \qquad \mathbf{E}\varepsilon = 0, \qquad \mathbf{E}\varepsilon^2 = \sigma^2, \qquad \mathbb{P}\{\varepsilon = b\} > 0.$$

It is easy to check that

$$\mathbb{P}\{\varepsilon = -\sigma^2/b\} = b^2/(b^2 + \sigma^2), \qquad \mathbb{P}\{\varepsilon = b\} = \sigma^2/(b^2 + \sigma^2).$$

Assuming (one-sided) boundedness of the differences $X_k$, in this article it is shown that up to an absolute constant factor the tail probability $\mathbb{P}\{M_n \geq x\}$ is dominated by the probability $\mathbb{P}\{S_n \geq x\}$. The result can be interpreted by saying that the behavior of tail probabilities of martingales is controlled in a very precise way by the simplest possible stochastic experiment—a series of eventually asymmetric coin tosses. This is not unexpected due to a common belief that Bernoulli random variables are those that are the most stochastic. It is less unexpected that one can provide a relatively simple proof of this fact.

For differences $X_k$ of a martingale $M_n$, consider the following boundedness condition: There exists a positive nonrandom $b > 0$ such that

$$(1.4) \qquad\qquad \mathbb{P}\{X_k \leq b\} = 1 \qquad \text{for } k = 1, \ldots, n.$$

For the conditional variances $s_k^2 = \mathbf{E}(X_k^2 | \mathcal{F}_{k-1})$ of differences $X_k$ of $M_n$, consider the following boundedness condition: There exist nonrandom $\sigma_k^2 \geq 0$ such that

$$(1.5) \qquad\qquad \mathbb{P}\{s_k^2 \leq \sigma_k^2\} = 1 \qquad \text{for } k = 1, \ldots, n.$$

THEOREM 1.1. *Assume that the differences $X_k$ of a martingale $M_n$ satisfy the conditions (1.4) and (1.5). Then, for all $x \in \mathbb{R}$, we have*

$$(1.6) \qquad\qquad \mathbb{P}\{M_n \geq x\} \leq \frac{e^2}{2} \mathbb{P}^{\circ}\{S_n \geq x\}$$

*with $e^2/2 \leq 3.7$, where $S_n$ is a sum of $n$ independent copies of a Bernoulli random variable $\varepsilon = \varepsilon(\sigma^2, b)$ with $\sigma^2 = (\sigma_1^2 + \cdots + \sigma_n^2)/n$ (the meaning of $\mathbb{P}^{\circ}$ is explained below). The inequality (1.6) yields*

$$(1.7) \qquad\qquad \mathbb{P}\{M_n \geq x\} \leq \frac{e^2}{2} \mathbb{P}^{\circ}\left\{\eta \geq \lambda + \frac{x}{b}\right\},$$

*where $\eta$ is a Poisson random variable with the parameter $\lambda = (\sigma_1^2 + \cdots + \sigma_n^2)/b^2$.*

The bound (1.7) is much rougher compared with (1.6) because it has to cover the case $n = \infty$, which supplies the heaviest tails. In general, tails of a sum of independent eventually nonidentically distributed Bernoulli random variables can have a complicated structure.

Let me explain the meaning of $\mathbb{P}^{\circ}$. Write $B(x) = \mathbb{P}\{S_n \geq x\}$ for the survival function of $S_n$. For $x$ such that $B(x) = 1$ or $B(x) = 0$ or when the function $B$ has a positive jump down, understand $\mathbb{P}^{\circ}$ just as probability. Let $B^{\circ}$ be a log-concave hull of $B$, that is, a minimal function such that $B \leq B^{\circ}$ and the function $x \mapsto -\log B^{\circ}(x)$ is a convex function. Define $\mathbb{P}^{\circ}\{S_n \geq x\} = B^{\circ}(x)$. It is easy to see (cf. Lemma 4.1) that in the case of the binomial or Poisson survival function $B$, the function $B^{\circ}$ is a log-linear interpolation of $B$: if $x < z < y$ and $x$ and $y$ are adjacent points where $B$ has positive jumps down, then

$$(1.8) \qquad B^{\circ}(z) = B^{1-\lambda}(x)\, B^{\lambda}(y), \qquad \text{if } z = (1-\lambda)x + \lambda y, \ 0 < \lambda < 1.$$

Similarly I introduce the linear interpolation $B^\diamond$ of $B$ by writing $B^\diamond(x) = B(x)$, for $x$ such that $B(x) = 1$ or $B(x) = 0$ or the function $B$ has a positive jump down, and

$$B^\diamond(z) = (1 - \lambda)B(x) + \lambda B(y) \qquad \text{for } x, y, z, \lambda \text{ as in } (1.8).$$

We have $B \leq B^\circ \leq B^\diamond$.

For differences $X_k$ of a martingale $M_n$, consider the boundedness condition

$$(1.9) \qquad \mathbb{P}\{-p_k \leq X_k \leq 1 - p_k\} = 1 \qquad \text{for } k = 1, \ldots, n,$$

where $p_k$ are nonrandom (it is clear that $0 \leq p_k \leq 1$).

THEOREM 1.2.   *Assume that the differences $X_k$ of a martingale $M_n$ satisfy the condition* (1.9). *Then, for $x \in \mathbb{R}$, we have*

$$(1.10) \qquad \mathbb{P}\{M_n \geq x\} \leq e\mathbb{P}^\circ\{S_n \geq x\}$$

*with $e \leq 2.72$, where $S_n = \varepsilon_1 + \cdots + \varepsilon_n$ is a sum of $n$ independent copies of a Bernoulli random variable*

$$\varepsilon = \varepsilon(p - p^2, 1 - p) \qquad \text{with } p = (p_1 + \cdots + p_n)/n.$$

*Furthermore, we have*

$$(1.11) \qquad \mathbb{P}\{M_n \geq x\} \leq \frac{e^3}{2}\mathbb{P}^\circ\left\{\eta \geq \lambda + \frac{x}{1 - p}\right\}$$

*with $e^3/2 \leq 10.1$, where $\eta$ is a Poisson random variable with $\lambda = pn/(1 - p)$.*

It is easy to check that the Bernoulli random variable $\varepsilon$ from Theorem 1.2 satisfies

$$\mathbb{P}\{\varepsilon = -p\} = 1 - p, \qquad \mathbb{P}\{\varepsilon = 1 - p\} = p.$$

By an application of (1.10) to $-M_n$, one can derive bounds for $\mathbb{P}\{M_n \leq x\}$.

For differences $X_k$ of a martingale $M_n$, consider the following boundedness conditions: There exist nonrandom $b_k \geq 0$ such that, for $k = 1, \ldots, n$,

$$(1.12) \qquad \mathbb{P}\{X_k \leq b_k\} = 1$$

and

$$(1.13) \qquad \mathbb{P}\{|X_k| \leq b_k\} = 1.$$

Write

$$(1.14) \qquad a_k = \max\{b_k, \sigma_k\}, \qquad a^2 = (a_1^2 + \cdots + a_n^2)/n,$$

where $\sigma_k$ are from the condition (1.5).

THEOREM 1.3. *Assume that the differences $X_k$ of a martingale $M_n$ satisfy the condition* (1.5) *and the one-sided boundedness condition* (1.12). *Then, for all $x \in \mathbb{R}$, we have*

$$(1.15) \qquad \mathbb{P}\{M_n \geq x\} \leq \frac{2e^3}{9} \mathbb{P}^{\circ}\{S_n \geq x\}$$

*with $2e^3/9 \leq 4.47$, where $S_n$ is a sum of $n$ independent copies of a symmetric Bernoulli random variable $\varepsilon = \varepsilon(a^2, a)$ with $a^2$ defined by* (1.14). *The inequality* (1.15) *implies*

$$(1.16) \qquad \mathbb{P}\{M_n \geq x\} \leq \frac{2e^3}{9}\left(1 - \Phi\left(\frac{x}{a}\right)\right),$$

*where $\Phi$ is the standard normal distribution function.*

The symmetric Bernoulli random variable from Theorem 1.3 satisfies $\mathbb{P}\{\varepsilon = \pm a\} = \frac{1}{2}$.

The following corollary is less general compared to Theorem 1.3.

COROLLARY 1.4. *Assume that the differences $X_k$ of a martingale $M_n$ satisfy the symmetric boundedness condition* (1.13). *Then, for all $x \in \mathbb{R}$, the bounds* (1.15) *and* (1.16) *of Theorem 1.3 hold, replacing $a^2$ with $(b_1^2 + \cdots + b_n^2)/n$.*

Theorems 1.1–1.3 show that the martingale type dependence does not influence the bounds for tail probabilities much compared to the independent, the i.i.d. and even the i.i.d. Bernoulli cases.

Most probably, the values of constants in Theorems 1.1–1.3 are not optimal; the preferred intention herein was to simplify the proofs as far as possible. A more powerful method that can improve constants and the structure of the bounds was used by Bentkus (2001). A bound from Bentkus (2001) applies to the special case of Theorem 1.2 when $p_k = 1/2$, and is precise for integer $x$ (a bound which is precise for all $x$ is in preparation). A consequence is that constants in the bounds (1.6), (1.10) and (1.15) of Theorems 1.1–1.3 cannot be smaller than 2, and these constants, say $c$, have to satisfy

$$2 \leq c \leq 3.7, \qquad 2 \leq c \leq 2.72, \qquad 2 \leq c \leq 4.47,$$

respectively, which means that space for improvement is restricted. In the case of Theorem 1.1, the multiplicative factor of losses in (1.6) is at most 1.85. In contrast to the martingale dependence, finding precise values of these constants in the independent and i.i.d. cases is considered a very difficult mathematical problem. For a given $n$, let $c_n$ be the best possible constant in Theorem 1.1. An impression that the sequence $c_n$ is increasing as $n \to \infty$ and that $\lim_{n\to\infty} c_n = 2$ is supported by the fact that $c_1 = 1.555884$.

Another supporting heuristic argument comes from the analysis of constants in the Berry–Esseen bounds in cases $n = 1$ and $n = \infty$ in Bentkus and Kirsha (1989) and Bentkus (1994), where a similar picture was observed.

One cannot replace $\mathbb{P}^\diamond\{S_n \geq x\}$ in Theorems 1.1–1.3 with $\mathbb{P}\{S_n \geq x\}$. This taboo clearly follows from the results (approach) of Bentkus (2001). A truly simple proof is provided in Section 4 as Lemma 4.8.

Let us compare Hoeffding's (1963) inequalities with bounds of Theorems 1.1–1.3. By Theorem 1 in Hoeffding (1963), under the conditions and notation of Theorem 1.2, we have

$$(1.17) \qquad \mathbb{P}\{M_n \geq x\} \leq H^n(p + x/n; p),$$

where, for $0 \leq p \leq 1$,

$$(1.18) \qquad H(a; p) = \left(\frac{1-p}{1-a}\right)^{1-a} \left(\frac{p}{a}\right)^a \qquad \text{for } p < a \leq 1$$

and

$$H(a; p) = 1 \qquad \text{for } a \leq p; \qquad H(a; p) = 0 \qquad \text{for } a > 1.$$

Among all inequalities that have a product form, Hoeffding's bounds are the best possible; see Lemma 4.7. Naturally, the product structure in our bounds is lost.

Hoeffding's inequalities remained unimproved until 1995 when Talagrand (1995a, b) inserted certain missing factors. Assuming independence and under the conditions of Theorem 1.2, Talagrand's bound is as follows: There exists an absolute constant $c > 0$ such that

$$(1.19) \quad \mathbb{P}\{M_n \geq x\} \leq \left(\frac{c\delta}{\delta + x} + \frac{c}{\delta}\right) H^n\left(p + \frac{x}{n}; p\right) \qquad \text{for } c \leq x \leq \frac{\delta^2}{c},$$

where $\delta^2 = np(1 - p)$. The right-hand side of (1.19) is simplified up to an absolute factor. This is a nonessential loss because Talagrand's bound depends on an inexplicit absolute constant.

I have a feeling that more or less explicit analytical functions do not truly follow the behavior of $\mathbb{P}\{M_n \geq x\}$ correctly.

The loss in Theorem 1.2 is at most the factor $e/2 \leq 1.36$. Up to an absolute constant, Theorem 1.2 (and Theorems 1.1 and 1.3 as well) says that the tail probability is maximized in the case of the simplest possible stochastic model, namely, in the case of a series of eventually asymmetric Bernoulli trials. To estimate $\mathbb{P}\{S_n \geq x\}$ one can use *special* bounds for the binomial probabilities [see, e.g., Shorack and Wellner (1986)], and in the view of Theorems 1.1–1.3, these special bounds are not so special at all.

Let us move on to Hoeffding's Theorem 3. To simplify notation (and without loss of generality) we assume that the number $b$ in Theorem 1.1 satisfies

$b = 1$. Assuming independence and under the conditions and notation of Theorem 1.1, Hoeffding proved that

$$(1.20) \qquad \mathbb{P}\{M_n \geq x\} \leq H^n\left(\frac{\sigma^2 + x/n}{1 + \sigma^2}; \frac{\sigma^2}{1 + \sigma^2}\right).$$

The simplest Hoeffding bound (1.17) is implied by (1.20) by using rescaling and choosing the maximal possible variance $\sigma^2 = p - p^2$ for distributions supported by the interval $[-p, 1-p]$. Assuming, in addition, that $|X_k| \leq B$, Talagrand (1995b) improved (1.20): There exists an absolute constant $c > 0$ such that

$$(1.21) \quad \mathbb{P}\{M_n \geq x\} \leq \left(\frac{c\sigma\sqrt{n}}{\sigma\sqrt{n} + x} + \frac{cB}{\sigma\sqrt{n}}\right) H^n\left(\frac{\sigma^2 + x/n}{1 + \sigma^2}; \frac{\sigma^2}{1 + \sigma^2}\right)$$

for $0 \leq x \leq n\sigma^2/(cB)$. Talagrand noticed that it is unclear how to improve (1.20) without assumptions like $|X_k| \leq B$. The inequality (1.21) nicely improves (1.20) when the variance is not too small, that is, in cases of Gaussian type behavior. To see this better, assume for simplicity that $B = 1$. Then, in the case of $\sigma^2 = 1$, the factor in (1.21) is on the order of $\sim \sqrt{n}/x$, in the range $\sqrt{n} \ll x \ll n$. However, for degenerating $\sigma^2 \to 0$ (i.e., when the behavior is of Poisson type), the range starts to shrink. To be definite, take $\sigma^2 = 1/n$. Then the factor is $\sim 1$ and $x$ has to satisfy $x \ll 1$. Notice that in such cases Theorem 1.1 still provides nice upper bounds in the whole range $x \leq n$ of interest.

Theorem 1.3 extends and refines Hoeffding's Theorem 2. The bound (1.16) with a somewhat worse constant is contained in Bentkus (2003). Using another approach, in Bentkus (2004) a bound similar to (1.16) was proved under the asymmetric boundedness condition

$$(1.22) \qquad \mathbb{P}\{d_k - a_k \leq X_k \leq d_k + a_k\} = 1,$$

where $d_k = d_k(X_1, \ldots, X_{k-1})$ are arbitrary $\mathcal{F}_{k-1}$-measurable random variables. This bound applies to the measure concentration. It is unclear whether one can extend and refine Hoeffding's Theorem 2 under the condition (1.22) using the methods of this article. Pinelis (1998) proved (1.16) under the symmetric boundedness condition (1.13). Earlier [see Pinelis (1999), Theorem 5], the bound (1.15) under the symmetric boundedness condition of Corollary 1.4 was established by Pinelis, assuming independence, for integer $x$ such that $x \in n + 2\mathbb{Z}$ and $|x| \leq n$.

Hoeffding's Theorem 2 had a considerable impact on research related to the measure concentration phenomena. For an introduction to the topic, see Gromov and Milman (1983), Alon and Milman (1984), Milman (1985, 1988), Milman and Schechtman (1986), McDiarmid (1989), Talagrand (1995a) and Ledoux (1999).

For statistical applications, optimal bounds for finite (i.e., fixed) $n$ are of interest [see Bentkus and van Zuijlen (2003)]. In this sense, the results herein are not optimal and hopefully can be improved by extending the methods of Bentkus (2001, 2004, b). However, the extensions involve considerable technical difficulties.

The history of inequalities for tail probabilities is a very rich classical topic [see, e.g., books Petrov (1975) and Shorack and Wellner (1986)]. The names Chernoff, Bennett, Prokhorov and Hoeffding come to mind. For $x \geq$ constant, the bounds above refine all the classical bounds. Indeed, one can estimate the binomial probability $\mathbb{P}\{S_n \geq x\}$ using these bounds.

*On methods.* Hoeffding (1963) applied the Chebyshev inequality to replace an indicator function of an interval with an exponential function, which can be interpreted as a kind of Fourier–Laplace transform. The further Hoeffding proof is precise; hence such a method cannot be used to improve the bounds. Talagrand (1995b) started with the Esscher transform, which is related to exponential functions. The proof in the articles by Pinelis [following Eaton (1970, 1974)] starts similarly to that of Hoeffding, but he used the functions $x \mapsto \max\{0; (x-t)^p\}$ instead of exponentials, with some $t \in \mathbb{R}$ and $p \in \mathbb{Z}$. In this article, we start in the same way. A nice and short argument in the proof of Theorem 1.2, which allows derivation of the inequality (3.9) from the bound (3.5), is extracted from an article by Pinelis (1999); see Lemma 4.2 below. The argument reduces the proof of Theorem 1.2 to the verification of inequalities (3.2) and (3.3). The scheme of the proof of Theorems 1.1 and 1.3 is similar, replacing (3.2) and (3.3) by appropriate counterparts. It seems that methods used here do not allow improvement of the constants and the structure of the bounds. In the aforementioned articles the methods rely on induction on $n$. Potentially such induction based methods can provide optimal bounds and it seems as well that they are more robust against generalizations.

**2. Some supplements, improvements and extensions.** In this section I provide some well-known upper bounds for Poisson and normal survival functions, a bit more complicated versions of bounds of Theorems 1.1–1.3 and precise bounds in the case $n = 1$.

A standard rather rough upper bound for a Poisson survival function is

$$\mathbb{P}\{\eta \geq \lambda + x\} \leq \exp\{x - (x + \lambda)\log(1 + x/\lambda)\} \qquad \text{for } x \geq 0.$$

For larger $x$, an impression about Poisson tails can provide the following inequalities [see Proposition 3 in Paulauskas (2002)]: There exist absolute positive constants $c_1$ and $c_2$ such that

$$c_1 g(x) \leq \mathbb{P}\{\eta \geq \lambda + x\} \leq c_2 g(x) \qquad \text{for } x \geq \max\{\lambda - 1, 1\},$$

where

$$g(x) = (\lambda + x)^{-1/2}(1 + x/\lambda)^{\{\lambda + x\} - 1} \exp\{x - (x + \lambda)\log(1 + x/\lambda)\}$$

and where $\{\lambda + x\}$ is the fractional part of $\lambda + x$.

A commonly used upper bound for the standard normal tail is

$$1 - \Phi(x) \leq \varphi(x)/x, \qquad \varphi(x) = (2\pi)^{-1/2}\exp\{-x^2/2\}, \qquad x > 0.$$

Let us pass to the extensions of Theorems 1.1–1.3. The extensions are more convenient in applications because they do not require checking of log-concavity. Hence, one can manipulate the bounds, for example, by applying limit theorems, and check log-concavity at the final stage of the application. I provide as well a direct generalization and extension of Hoeffding's Theorem 3 to martingales. This extension can be useful in cases where checking log-concavity is not available or in cases when very precise bounds are not needed. It is interesting to notice that in contrast to the much more subtle and powerful Theorem 3, there exist lots of extensions, improvements and generalizations of Hoeffding's Theorem 2. Probably the reason is that Theorem 2 is simpler than Theorem 3, because instead of variances, it involves only rather rough size parameters.

For differences $X_k$ of a martingale $M_n$, consider the following boundedness condition: There exist positive nonrandom $b_k > 0$ and $\sigma_k > 0$ such that

$$(2.1) \qquad \mathbb{P}\{X_k \leq b_k\} = 1, \qquad \mathbb{P}\{s_k^2 \leq \sigma_k^2\} = 1$$

for $k = 1, \ldots, n$, where $s_k^2 = \mathbf{E}(X_k^2 | \mathcal{F}_{k-1})$ are the conditional variances of $X_k$.

Introduce independent Bernoulli random variables $\theta_k = \theta_k(\sigma_k^2, b_k)$ such that

$$(2.2) \qquad \mathbf{E}\theta_k = 0, \qquad \mathbf{E}\theta_k^2 = \sigma_k^2, \qquad \mathbb{P}\{\theta_k = b_k\} > 0$$

[cf. the definition (1.3) of the Bernoulli random variable $\varepsilon = \varepsilon(\sigma, b)$]. Write

$$(2.3) \qquad T_n = \theta_1 + \cdots + \theta_n.$$

THEOREM 2.1. *Assume that the differences $X_k$ of a martingale $M_n$ satisfy (2.1). Then, for $h > 0$, we have*

$$(2.4) \qquad \mathbb{P}\{M_n \geq x\} \leq \exp\{-hx\}\mathbf{E}\exp\{hT_n\}$$

*for all $x \in \mathbb{R}$. If all $b_k$ are equal, $b_k = b$, then*

$$(2.5) \quad \mathbb{P}\{M_n \geq x\} \leq \inf_{h > 0} \exp\{-hx\}\mathbf{E}\exp\{hS_n\} = H^n\left(\frac{\sigma^2 + bx/n}{b^2 + \sigma^2}; \frac{\sigma^2}{b^2 + \sigma^2}\right).$$

*Here $S_n$ is a sum of $n$ independent copies of a Bernoulli random variable $\varepsilon = \varepsilon(\sigma^2, b)$ with the variance $\sigma^2 = (\sigma_1^2 + \cdots + \sigma_n^2)/n$ and the function $H$ is given by (1.18).*

We provide proofs in Section 3. A number of upper bounds for the function $H$ are provided in Hoeffding (1963).

Let $x_+ = \max\{0, x\}$ and $x_+^s = (x_+)^s$.

THEOREM 2.2.   *Write* $f(x) = (x - t)_+^s$, *where* $s \geq 2$, *and assume that the differences* $X_k$ *of a martingale* $M_n$ *satisfy* (2.1). *Then, for all* $t < x$ *and* $x \in \mathbb{R}$, *we have*

$$(2.6) \qquad \mathbb{P}\{M_n \geq x\} \leq \mathbf{E} f(T_n)/(x - t)^s$$

*and*

$$(2.7) \qquad \mathbb{P}\{M_n \geq x\} \leq e^s s^{-s} \Gamma(s + 1) \mathbb{P}^\circ\{T_n \geq x\},$$

*where* $x \mapsto \mathbb{P}^\circ\{T_n \geq x\}$ *is a log-concave hull of the survival function* $x \mapsto \mathbb{P}\{T_n \geq x\}$ *and* $\Gamma$ *is the gamma function.*

In the next proposition we provide precise bounds for $n = 1$ under the conditions of Theorems 1.1–1.3.

PROPOSITION 2.3 (Case $n = 1$).   *Assume that a random variable satisfies* $\mathbf{E} X = 0$. *Let* $a < 0 < b$ *and* $\sigma > 0$.

(i) *Let* $B(x) = \sup \mathbb{P}\{X \geq x\}$, *where* sup *is taken over all random variables* $X$ *such that* $\mathbb{P}\{a \leq X \leq b\} = 1$. *Then* $B(x) = p$ *with* $p = -a/(x - a)$ *for* $0 \leq x \leq b$.

(ii) *Write* $B(x) = \sup \mathbb{P}\{X \geq x\}$, *where* sup *is taken over all* $X$ *such that*

$$\mathbb{P}\{X \leq b\} = 1 \quad and \quad \mathbf{E} X^2 \leq \sigma^2.$$

*Then* $B(x) = p$ *with* $p = \sigma^2/(x^2 + \sigma^2)$ *for* $0 \leq x \leq b$.

*In both cases* (i) *and* (ii) *we have* $B(x) = 1$ *for* $x \leq 0$ *and* $B(x) = 0$ *for* $x > b$.

**3. Proofs.** Write $x_+ = \max\{0, x\}$ and $x_+^s = (x_+)^s$. Let us start with the proof of Theorem 1.2 because it is simpler compared to the proof of Theorem 1.1.

PROOF OF THEOREM 1.2.   Let us prove first the bound (1.10). In the proof we assume that $-pn < x \leq n - pn$ because for other values of $x$ the inequality (1.10) reduces either to $1 \leq e$ or to $0 \leq 0$, which is obvious.

Write $f(z) = (z - t)_+$, where $t \in \mathbb{R}$ is a parameter to be chosen later. Notice that $\mathbb{I}\{u \geq x\} \leq f(u)/(x - t)$ for $t < x$, where $\mathbb{I}\{A\}$ is the indicator function of event $A$. Using the Chebyshev inequality, we have

$$(3.1) \qquad \mathbb{P}\{M_n \geq x\} \leq \mathbf{E} f(M_n)/(x - t) \qquad \text{for } t < x.$$

Applying Lemma 4.3, we have

$$\mathbf{E}f(M_n) \leq \mathbf{E}f(T_n), \tag{3.2}$$

where $T_n = \xi_1 + \cdots + \xi_n$ is a sum of independent Bernoulli random variables $\xi_k$ such that

$$\mathbb{P}\{\xi_k = -p_k\} = 1 - p_k \quad \text{and} \quad \mathbb{P}\{\xi_k = 1 - p_k\} = p_k.$$

We are going to replace the eventually non-i.i.d. Bernoulli random variables $\xi_k$ with the i.i.d. Bernoulli random variables from the condition of the theorem. If $f$ is a convex function, then

$$\mathbf{E}f(T_n) \leq \mathbf{E}f(\varepsilon_1 + \cdots + \varepsilon_n) = \mathbf{E}f(S_n), \tag{3.3}$$

with $S_n$ from the condition of the theorem. Hoeffding [(1956), Theorem 3] proved (3.3) for strictly convex $f$ and Gleser [(1975), Corollary 2.1] extended (3.3) to convex $f$. One can easily check (3.3) using the Schur concavity; see the proof of Lemma 4.5 for a definition of Schur concave functions.

In the specific case of $f(x) = (x - t)_+$ we have

$$\mathbf{E}f(S_n) = -\int_t^\infty (z - t)\, d\mathbb{P}\{S_n \geq z\} = \int_t^\infty \mathbb{P}\{S_n \geq z\}\, dz. \tag{3.4}$$

Combining (3.1)–(3.4), we obtain

$$\mathbb{P}\{M_n \geq x\} \leq \inf_{t < x} \frac{1}{x - t} \int_t^\infty \mathbb{P}\{S_n \geq z\}\, dz. \tag{3.5}$$

To estimate the right-hand side of (3.5), we can apply Lemma 4.2 with

$$\alpha = -pn, \qquad \beta = n - pn, \qquad s = 1 \quad \text{and} \quad B(z) = \mathbb{P}\{S_n \geq z\}.$$

We get $\mathbb{P}\{M_n \geq x\} \leq e\mathbb{P}^\circ\{S_n \geq z\}$, which concludes the proof of (1.10).

Let us prove prove (1.11). The sum $S_n$ is a sum of $n$ independent copies of a Bernoulli random variable $\varepsilon = \varepsilon(p - p^2, 1 - p)$. Applying the bound (1.7) of Theorem 1.1, we obtain

$$\mathbb{P}\{S_n \geq z\} \leq \frac{e^2}{2} \mathbb{P}^\circ\left\{\eta \geq \lambda + \frac{z}{1 - p}\right\}, \tag{3.6}$$

where $\eta$ is a Poisson random variable with $\lambda = pn/(1 - p)$. Combining inequalities (3.5) and (3.6), we have

$$\mathbb{P}\{M_n \geq x\} \leq \frac{e^2}{2} \inf_{t < x} \frac{1}{x - t} \int_t^\infty \mathbb{P}^\circ\left\{\eta \geq \lambda + \frac{z}{1 - p}\right\} dz, \tag{3.7}$$

and an application of Lemma 4.2 yields (1.11). $\quad\square$

PROOF OF THEOREM 2.1. Let us prove the bound (2.4). Using the Chebyshev inequality, we have

$$\mathbb{P}\{M_n \geq x\} \leq \exp\{-hx\}\mathbf{E}\exp\{hM_n\} \qquad \text{for } h > 0.$$

By Lemma 4.4, we have $\mathbf{E}\exp\{hM_n\} \le \mathbf{E}\exp\{hT_n\}$, which concludes the proof of (2.4).

Let us prove (2.5). The inequality

$$\exp\{-hx\}\mathbf{E}\exp\{hT_n\} \le \exp\{-hx\}\mathbf{E}\exp\{hS_n\}$$

is proved in Hoeffding [(1963), (4.22) in the proof of Theorem 3]. The equality in (2.5) is just the definition of the Hoeffding function; see Hoeffding (1963). □

PROOF OF THEOREM 2.2. Let us prove (2.6). Write $f(z) = (z-t)_+^s$. Using the Chebyshev inequality, we have $P\{M_n \ge x\} \le \mathbf{E}f(M_n)/(x-t)^s$ for $t < x$. By Lemma 4.4, we can estimate $\mathbf{E}f(M_n) \le \mathbf{E}f(T_n)$, and (2.6) follows. The bound (2.7) is implied by Lemma 4.2. □

PROOF OF THEOREM 1.1. Let us prove (1.6). Without loss of generality (rescaling if necessary), we can assume that the number $b$ from the condition (1.4) satisfies $b = 1$. In the proof we assume that $-n\sigma^2 < x \le n$, because for $-n\sigma^2 < x$ or $x > n$ the inequality (1.6) reduces to obvious $1 \le e^2/2$ or $0 \le 0$, respectively.

Write $f(z) = (z-t)_+^2$, where $t \in \mathbb{R}$ is a parameter to be chosen later. Using the Chebyshev inequality, we have

$$(3.8) \qquad \mathbb{P}\{M_n \ge x\} \le \mathbf{E}f(M_n)/(x-t)^2 \qquad \text{for } t < x.$$

By Lemma 4.4, we can replace $M_n$ with a sum of Bernoulli random variables, that is,

$$(3.9) \qquad \mathbf{E}f(M_n) \le \mathbf{E}f(T_n),$$

where $T_n = \theta_1 + \cdots + \theta_n$ is a sum of independent eventually nonidentically distributed Bernoulli random variables $\theta_k = \theta(\sigma_k^2, 1)$ [cf. (2.2)].

By Lemma 4.5, we can replace the non-i.i.d. with the i.i.d. Bernoulli random variables,

$$(3.10) \qquad \mathbf{E}f(T_n) \le \mathbf{E}f(\varepsilon_1 + \cdots + \varepsilon_n) = \mathbf{E}f(S_n),$$

where $S_n = \varepsilon_1 + \cdots + \varepsilon_n$ is a sum of i.i.d. Bernoulli random variables $\varepsilon_k = \varepsilon_k(\sigma^2, 1)$.

In the specific case of $f(x) = (x-t)_+^2$, we have

$$(3.11)\, \mathbf{E}f(S_n) = -\int_t^\infty (z-t)^2 \, d\mathbb{P}\{S_n \ge z\} = 2\int_t^\infty (z-t)\mathbb{P}\{S_n \ge z\} \, dz.$$

Combining (3.8)–(3.11), we obtain

$$(3.12) \qquad \mathbb{P}\{M_n \ge x\} \le \inf_{t<x} \frac{2}{(x-t)^2} \int_t^\infty (z-t)\mathbb{P}\{S_n \ge z\} \, dz.$$

To estimate the right-hand side of (3.12), we can apply Lemma 4.2 with

$$\alpha = -n\sigma^2, \qquad \beta = n, \qquad s = 2 \quad \text{and} \quad B(z) = \mathbb{P}\{S_n \geq z\}.$$

We get $\mathbb{P}\{M_n \geq x\} \leq (e^2/2)\mathbb{P}^\circ\{S_n \geq z\}$, proving (1.6).

It remains to prove (1.7). Introduce the martingale $K_{n+m} = Y_1 + \cdots + Y_{n+m}$ with the differences

$$
\begin{aligned}
(3.13) \qquad & Y_k = X_k, \qquad \text{for } k = 1, \ldots, n \quad \text{and} \\
& Y_k = 0, \qquad \text{for } k = n+1, \ldots, n+m.
\end{aligned}
$$

To the martingale $K_{n+m}$ we can apply the bound (1.6) of Theorem 1.1. We get

$$(3.14) \qquad \mathbb{P}\{M_n \geq x\} = \mathbb{P}\{K_{n+m} \geq x\} \leq \frac{e^2}{2}\mathbb{P}^\circ\{S_{n+m} \geq x\},$$

where $S_{n+m}$ is a sum of $n+m$ independent copies of a Bernoulli random variable

$$\varepsilon = \varepsilon(\sigma^2, b) \qquad \text{with } \sigma^2 = (\sigma_1^2 + \cdots + \sigma_n^2)/(n+m) = b^2\lambda/(n+m).$$

Centering and rescaling, we get

$$(3.15) \qquad \mathbb{P}\{S_{n+m} \geq x\} = \mathbb{P}\{Z_{n+m} \geq (\lambda + x/b)/(1 + \lambda/(n+m))\},$$

where $Z_{n+m}$ is a sum of $n+m$ independent copies of a Bernoulli random variable, say $\xi$, such that

$$\mathbb{P}\{\xi = 0\} = q \qquad \text{with } q = 1 - p \quad \text{and} \quad p = \mathbb{P}\{\xi = 1\} = \lambda/(n+m+\lambda).$$

To the sum $Z_{n+m}$ we can apply the Poisson limit theorem because $p(n+m) \to \lambda$ as $m \to \infty$. We get

$$(3.16) \quad \lim_{m\to\infty} \mathbb{P}\{Z_{n+m} \geq (\lambda + x/b)/(1 + \lambda/(n+m))\} = \mathbb{P}\{\eta \geq \lambda + x/b\}.$$

Combining (3.14)–(3.16), we conclude the proof of (1.7). $\square$

PROOF OF THEOREM 1.3. Let us prove (1.15). Write $f(z) = (z - t)_+^3$. Similar to (3.1) and (3.8), we have

$$(3.17) \qquad \mathbb{P}\{M_n \geq x\} \leq \mathbf{E}f(M_n)/(x-t)^3 \qquad \text{for } t < x.$$

Let $T_n = \theta_1 + \cdots + \theta_n$ be a sum of independent Bernoulli random variables such that $\mathbb{P}\{\theta_k = -a_k\} = \mathbb{P}\{\theta_k = a_k\} = 1/2$. An application of Lemma 4.6 yields $\mathbf{E}f(M_n) \leq \mathbf{E}f(T_n)$. The inequality

$$(3.18) \qquad \mathbf{E}f(T_n) \leq \mathbf{E}f(\varepsilon_1 + \cdots + \varepsilon_n) = \mathbf{E}f(S_n),$$

where $S_n = \varepsilon_1 + \cdots + \varepsilon_n$ is a sum of $n$ independent copies of a symmetric Bernoulli random variable $\varepsilon = \varepsilon(a^2, a)$ as established in Eaton ([1970], [1974]) and Pinelis ([1994]).

Using $f(x) = (x - t)^3$, integrating by parts and combining (3.17) and (3.18), we have

$$\mathbb{P}\{M_n \geq x\} \leq \inf_{t < x} \frac{3}{(x - t)^3} \int_t^\infty (z - t)^2 \mathbb{P}\{S_n \geq z\} \, dz.$$

Now an application of Lemma 4.2 implies (1.15).

It remains to prove (1.16). Introduce the martingale $K_{n+m} = Y_1 + \cdots + Y_{n+m}$ with the differences defined by (3.13). To the martingale $K_{n+m}$ we can apply the bound (1.15) of Theorem 1.3. We get

$$(3.19) \quad \mathbb{P}\{M_n \geq x\} = \mathbb{P}\{K_{n+m} \geq x\} \leq \frac{2e^3}{9} \mathbb{P}^\circ \left\{ (n + m)^{-1/2} S_{n+m} \geq \frac{x}{a} \right\},$$

where $S_{n+m}$ is a sum of $n + m$ independent copies of a symmetric Bernoulli random variable, say $\varepsilon$, such that $\mathbb{P}\{\varepsilon = -1\} = \mathbb{P}\{\varepsilon = 1\} = 1/2$. We conclude the proof of (1.16) by passing to the limit in (3.19) as $m \to \infty$ and using the central limit theorem. $\quad\square$

PROOF OF COROLLARY 1.4. The boundedness condition (1.13) guarantees that the conditional variances $s_k^2$ are bounded from above by $b_k^2$. Hence, $a_k = b_k$ and we can apply (1.15) and (1.16) with $a^2 = b_1^2 + \cdots + b_n^2$. $\quad\square$

PROOF OF PROPOSITION 2.3. It suffices to prove (i) and (ii) only for $0 < x \leq b$. Indeed, for $x > b$ we have obviously $B(x) = 0$. For $x \leq 0$, the upper bound $B(x) \leq 1$ is obvious; the lower bound $B(x) \geq 1$ follows by considering the random variable $X = 0$.

(i) For $0 \leq x \leq b$, the linear function $u(t) = (1 - p)t + p$ satisfies $\mathbb{I}\{t \geq x\} \leq u(t)$ for all $t$ from the interval $[a, b]$. Therefore, we have

$$(3.20) \qquad \mathbb{P}\{X \geq x\} \leq \mathbb{E}u(X) = p \quad \text{and} \quad B(x) \leq p.$$

The lower bound $B(x) \geq p$ is realized by a Bernoulli random variable, say $X = \varepsilon$, such that $\mathbb{P}\{\varepsilon = a\} = 1 - p$ and $\mathbb{P}\{\varepsilon = b\} = p$.

(ii) For $0 < x \leq b$ and $t \leq b$, the quadratic function $u(t) = (1 - p)^2 x^2 (t + \sigma^2/x)^2$ satisfies the inequality $\mathbb{I}\{t \geq x\} \leq u(t)$. Similar to (3.20) it follows that $B(x) \leq p$. The lower bound $B(x) \geq p$ is realized by a Bernoulli random variable, say $X = \varepsilon$, such that $\mathbb{P}\{\varepsilon = -\sigma^2/b\} = 1 - p$ and $\mathbb{P}\{\varepsilon = b\} = p$. $\quad\square$

**4. Auxiliary results.** A function $f : A \to [0, \infty)$ defined on a subset $A \subset R$ is called log-concave if the function $x \mapsto -\log f(x)$ is convex. Whereas $f$ can assume the value 0, the function $-\log f$ can assume the value $\infty$. Call a random variable $X$ discrete if there exists a countable set $A$ such that

$$\mathbb{P}\{X \in A\} = 1 \quad \text{and} \quad \mathbb{P}\{X = x\} > 0 \qquad \text{for all } x \in A.$$

A survival function $x \mapsto \mathbb{P}\{X \geq x\}$ is called discrete if $X$ is discrete. A binomial survival function $x \mapsto \mathbb{P}\{S_n \geq x\}$ is discrete and is not log-concave as a function defined on $\mathbb{R}$. However, it is log-concave as a function defined on the set $A$ of points at which it has positive jumps down (see Lemma 4.1). Therefore, a discrete survival function is called log-concave if it is log-concave as a function defined on the set $A$ (hopefully this terminology will not lead to misunderstanding). For a function $f : \mathbb{R} \to [0, \infty)$, introduce its log-concave hull $f^\circ : \mathbb{R} \to [0, \infty)$ as a minimal log-concave function such that $f \leq f^\circ$. Any survival function has a unique log-concave hull which is again a log-concave survival function.

For a random variable $X$, which assumes integer values, the probability mass function is defined as $p_n = \mathbb{P}\{X = n\}$ for $n \in \mathbb{Z}$. In the literature, distributions with log-concave densities and probability mass functions are refered to as strong unimodal in the sense of Ibragimov [cf. Keilson and Gerber (1971) and Ibragimov (1956)]. We are interested in log-concave survival functions, which have a weaker requirement compared to the strong unimodality. The next lemma is just a reexposition of some facts from Keilson and Gerber (1971) and Pinelis (1998, 1999).

LEMMA 4.1. (i) *Let $n \mapsto p_n$ and $n \mapsto q_n$ be log-concave functions such that $p_n, q_n \geq 0$. Then the convolution*

$$(p * q)_n = \sum_{k=-\infty}^{\infty} p_{n-k} q_k$$

*is a log-concave function.*

(ii) *Let $n \mapsto p_n$ be a log-concave function such that $p_n \geq 0$. Then the function $n \mapsto t_n$ with $t_n = \sum_{k \geq n} p_k$ is a log-concave function.*

(iii) *Bernoulli random variables have log-concave probability mass functions. Binomial survival functions are log-concave (as discrete ones).*

(iv) *Let $B_k$ be a sequence of log-concave survival functions which have probability mass functions supported by $\mathbb{Z}$. Then the pointwise limit $\lim_{k \to \infty} B_k$ is a log-concave function.*

(v) *Poisson survival functions are log-concave (as discrete ones).*

(vi) *Binomial and Poisson survival functions $B$ satisfy $B \leq B^\circ \leq B^\diamond$. In both cases $B^\circ$ is just a log-linear interpolation of $B$.*

In general, it is not true that a sum of two independent discrete random variables with log-concave survival functions has a log-concave discrete survival function. Indeed, let $\varepsilon, \varepsilon_1, \varepsilon_2$ be i.i.d. Bernoulli random variables such that

$$\mathbb{P}\{\varepsilon = 0\} = q \quad \text{and} \quad \mathbb{P}\{\varepsilon = 1\} = p \qquad \text{with } p + q = 1.$$

Then the discrete survival function $B(x) = \mathbb{P}\{\varepsilon_1 + a\varepsilon_2 \geq x\}$ is not log-concave provided that the numbers $a > 0$ and $p > 0$ are sufficiently small. Indeed, assume that the function $B$ is log-concave. The random variable $\varepsilon_1 + a\varepsilon_2$ assumes values $0 < a < 1 < 1 + a$. The log-concavity of $B$ yields $B(0)^{1-a}B(1)^a \leq B(a)$, which is equivalent to $p^a \leq 2p - p^2$. Passing to the limit as $a \downarrow 0$, we have $1 \leq 2p - p^2$, which is impossible if $p > 0$ is sufficiently small. A similar consideration shows that survival functions of discrete infinite divisible random variables are not necessarily log-concave: for example, the survival function of $\eta + a\xi$, where $\eta$ and $\xi$ are Poisson random variables with parameters $\lambda > 0$ and $\gamma > 0$ is not log-concave provided that $a > 0$ and $\gamma > 0$ are sufficiently small (just consider the values of the survival function at points $0, a$ and $2a$).

PROOF OF LEMMA 4.1. (i) Write

$$\delta = (p * q)_n^2 - (p * q)_{n-1}(p * q)_{n+1}.$$

We have to prove that $\delta \geq 0$. It is easy to check that $2\delta = \sum_{k,r=-\infty}^{\infty} \alpha\beta$ with

$$\alpha = p_k p_r - p_{k+1} p_{r-1} \quad \text{and} \quad \beta = q_{n-k} q_{n-r} - q_{n-k-1} q_{n-r+1}.$$

If $k \geq r$, then $\alpha \geq 0$ and $\beta \geq 0$, because both functions $n \mapsto p_n$ and $n \mapsto q_n$ are log-concave. If $k < r$, then $\alpha \leq 0$ and $\beta \leq 0$, which concludes the proof of $\delta \geq 0$.

(ii) Notice that $t_n = (p * q)_n$, where $q_n = \mathbb{I}\{n \leq 0\}$ is log-concave function, and apply (i).

(iii) It is clear that Bernoulli probability mass functions are log-concave. Hence, by applying (i), binomial probability mass functions are log-concave. Therefore, (ii) guarantees that binomial survival functions are log-concave.

(iv) Obvious.

(v) A Poisson survival function is a limit of a sequence of Bernoulli survival functions. Therefore we can apply (iii) and (iv).

(vi) The inequality $B \leq B^\circ$ is obvious. The inequality $B^\circ \leq B^\diamond$ is equivalent to the elementary inequality

$$a^{1-\lambda}b^\lambda \leq a + b \qquad \text{for } a \geq b \geq 0 \quad \text{and} \quad 0 \leq \lambda \leq 1.$$

To see that $\mathbb{B}^\circ$ is a log-linear interpolation of $B$, it suffices to compare the graphs of $-\log B$ and $-\log B^\circ$. $\square$

In the case of a log-concave $B$, the next lemma was proved by Pinelis (1999). Actually, the work by Pinelis contains more general results. In the special case $s = 1$, the result was established by Bretagnolle (1980) and by Kemperman [see Shorack and Wellner (1986), Chapter 25, Lemma 1].

LEMMA 4.2.  *Let $s > 0$. Let $B$ be a survival function with a log-concave hull $B^\circ$. Let*

$$\alpha = \sup\{y : B^\circ(y) = 1\} \quad and \quad \beta = \inf\{y : B^\circ(y) = 0\}.$$

*Then, for $x$ such that $\alpha < x \le \beta$, we have*

$$(4.1) \qquad \inf_{t < x}(x - t)^{-s} \int_t^\infty s(z - t)^{s-1} B(z)\, dz \le e^s s^{-s} \Gamma(s+1) B^\circ(x),$$

*where $\Gamma(s) = \int_0^\infty \tau^{s-1} \exp\{-\tau\}\, d\tau$.*

PROOF.  Because it is short, the proof is provided. The function $z \mapsto -\log \mathbb{B}^\circ(z)$ is a convex function. It is clear that this function is strictly positive and strictly increasing in the interval $(\alpha, \beta]$. Hence, for each $x \in (\alpha, \beta]$, there exists a linear function, say $y(z) = a + bz$, with some positive $b > 0$, such that

$$y(x) = -\log \mathbb{B}^\circ(x) \quad \text{and} \quad -\log \mathbb{B}^\circ(z) \ge y(z) \qquad \text{for all } z \in \mathbb{R}.$$

The numbers $a = a(x, B)$ and $b = b(x, B)$ can depend on $x$ and $B$. In particular, we have

$$(4.2) \quad \mathbb{B}^\circ(x) = \exp\{-a - bx\}, \qquad \mathbb{B}^\circ(z) \le \exp\{-a - bz\} \qquad \text{for all } z \in \mathbb{R}.$$

Using $B \le B^\circ$ and (4.2), we have

$$
\begin{aligned}
\int_t^\infty s(z - t)^{s-1} B(z)\, dz &\le \exp\{-a\} \int_t^\infty s(z - t)^{s-1} \exp\{-bz\}\, dz \\
(4.3) \qquad\qquad &= \Gamma(s+1) b^{-s} \exp\{-a - bt\} \\
&= \Gamma(s+1) b^{-s} \exp\{b(x - t)\} \mathbb{B}^\circ(x).
\end{aligned}
$$

Using (4.3) and choosing $t$ such that $b(x - t) = s$, we obtain (4.1).  □

It seems that the next lemma has to be a well-known fact [a useful related reference is Karlin and Studden (1966)]. We write $\xi = \xi(a, b)$ for a Bernoulli random variable such that

$$(4.4) \qquad \mathbb{P}\{\xi = a\} = b/(b - a) \quad \text{and} \quad \mathbb{P}\{\xi = b\} = -a/(b - a).$$

LEMMA 4.3.  (i) *Let $f: \mathbb{R} \to \mathbb{R}$ be a convex function. Assume that a random variable $X$ satisfies*

$$\mathbf{E}X = 0, \qquad \mathbb{P}\{a \leq X \leq b\} = 1, \qquad a \leq 0 \leq b.$$

*Then $\mathbf{E}f(X) \leq \mathbf{E}f(\xi)$, where $\xi$ is a Bernoulli random variable satisfying* (4.4).

(ii) *Let a function $f: \mathbb{R}^n \to \mathbb{R}$ be a convex function of each of variables $x_1, \ldots, x_n$ when the remaining $n - 1$ variables are kept fixed. Assume that the differences $X_k$ of a martingale $M_n = X_1 + \cdots + X_n$ satisfy*

$$\mathbb{P}\{a_k \leq X_k \leq b_k\} = 1,$$

*where numbers $a_k \leq 0 \leq b_k$ are nonrandom for all $k$. Let $\xi_k = \xi_k(a_k, b_k)$ be independent Bernoulli random variables. Then we have*

(4.5)                    $$\mathbf{E}f(X_1, \ldots, X_n) \leq \mathbf{E}f(\xi_1, \ldots, \xi_n).$$

PROOF.  (i) We have to prove that $\mathbf{E}f(X) \leq \mathbf{E}f(\xi)$. Let $u: [a, b] \to \mathbb{R}$ be a linear function. Then $\mathbf{E}u(X) = \mathbf{E}u(\xi)$ because $\mathbf{E}X = \mathbf{E}\xi = 0$. Choose $u$ such that $u(a) = f(a)$ and $u(b) = f(b)$. Then $f \leq u$ because $f$ is convex. Futhermore, $\mathbf{E}u(\xi) = \mathbf{E}f(\xi)$ because $\mathbb{P}\{\xi \in \{a, b\}\} = 1$. Hence $\mathbf{E}f(X) \leq \mathbf{E}u(X) = \mathbf{E}u(\xi) = \mathbf{E}f(\xi)$, which concludes the proof in the case (i).

(ii) We use induction in $n$. In the case $n = 1$, the result was proved in (i). Let $n > 1$ and let (4.5) hold for $1, \ldots, n - 1$. Notice that for given $X_1$, the sequence

(4.6)              $$Z_0 = 0, \qquad Z_1 = X_2, \ldots, Z_{n-1} = X_2 + \cdots + X_n$$

is a martingale sequence with differences that satisfy

(4.7)        $$\mathbb{P}\{a_{k+1} \leq Z_k - Z_{k-1} \leq b_{k+1}\} = 1 \qquad \text{for } k = 1, \ldots, n - 1.$$

Conditioning on $X_1$ and applying the induction assumption twice (for $n - 1$ and 1), we have

$$\begin{aligned}
\mathbf{E}f(X_1, \ldots, X_n) &= \mathbf{E}(f(X_1, \ldots, X_n)|X_1) \\
&\leq \mathbf{E}(f(X_1, \xi_2, \ldots, \xi_n)|X_1) \\
&= \mathbf{E}(f(X_1, \xi_2, \ldots, \xi_n)|\xi_2, \ldots, \xi_n) \\
&\leq \mathbf{E}(f(\xi_1, \xi_2, \ldots, \xi_n)|\xi_2, \ldots, \xi_n) \\
&= \mathbf{E}f(\xi_1, \ldots, \xi_n),
\end{aligned}$$

which completes the proof of (4.5) for $n > 1$.  $\square$

LEMMA 4.4. *Let $f$ be one of the functions*

$$f(x) = (x-t)_+^2, \qquad t \in \mathbb{R};$$
$$f(x) = (x-t)_+^s, \qquad s > 2;$$
$$f(x) = \exp\{hx\}, \qquad h > 0.$$

(i) *Assume that a random variable $X$ satisfies*

$$\mathbb{P}\{X \leq b\} = 1, \qquad \mathbf{E}X^2 \leq \sigma^2.$$

*Then $\mathbf{E}f(X) \leq \mathbf{E}f(\theta)$, where a Bernoulli random variable $\theta$ satisfies $\theta = \theta(\sigma^2, b)$ [see the definition (2.2) of $\theta$].*

(ii) *Assume that a martingale $M_n$ satisfies condition (2.1), that is, that $X_k \leq b_k$ and $s_k^2 \leq \sigma_k^2$ with probability 1, where $s_k^2$ are the conditional variances of the differences $X_k$. Let $T_n = \theta_1 + \cdots + \theta_n$ be a sum of independent Bernoulli random variables $\theta_k = \theta_k(\sigma_k^2, b_k)$. Then we have $\mathbf{E}f(M_n) \leq \mathbf{E}f(T_n)$.*

PROOF. It suffices to prove the lemma with $f(x) = (x-t)_+^2$, $t \in \mathbb{R}$. Indeed, both functions $g(x) = (x-t)_+^s$ and $g(x) = \exp\{hx\}$ with $s > 2$ and $h > 0$ allow the integral representation

$$(4.8) \qquad g(x) = \tfrac{1}{2} \int_{\mathbb{R}} g'''(u)(x-u)_+^2 \, du, \qquad g''' \geq 0.$$

Therefore, the inequality $\mathbf{E}(M_n - u)_+^2 \leq \mathbf{E}(T_n - u)_+^2$ for all $u \in \mathbb{R}$ clearly implies $\mathbf{E}g(M_n) \leq \mathbf{E}g(T_n)$.

Henceforth let $f(x) = (x-t)_+^2$.

(i) Let us prove that $\mathbf{E}f(X) \leq \mathbf{E}f(\theta)$. The r.v. $X$ satisfies $\mathbb{P}\{X \leq b\} = 1$. We consider the following cases separately:

(a) $t \leq -\sigma^2/b$;
(b) $-\sigma^2/b < t < b$;
(c) $t \geq b$.

*Case* (a). Using

$$(x-t)_+^2 \leq (x-t)^2 \quad \text{and} \quad \mathbf{E}X = 0, \qquad \mathbf{E}X^2 \leq \sigma^2,$$

we have

$$\mathbf{E}f(X) \leq \mathbf{E}(X-t)^2 \leq \sigma^2 + t^2 = \mathbf{E}(\theta-t)^2 = \mathbf{E}(\theta-t)_+^2 = \mathbf{E}f(\theta).$$

*Case* (b). Notice that

$$(x-t)_+^2 \le c(x+\sigma^2/b)^2, \qquad \text{for } x \le b, \text{ where } c = b^2(b-t)^2/(b^2+\sigma^2)^2.$$

Using this inequality and $\mathbf{E}X = 0$, $\mathbf{E}X^2 \le \sigma^2$, we obtain

$$\mathbf{E}f(X) \le c\mathbf{E}(X+\sigma^2)^2 \le c(\sigma^2 + \sigma^4/b^2) = \mathbf{E}(\theta-t)_+^2 = \mathbf{E}f(\theta).$$

*Case* (c). Now we have $\mathbf{E}f(X) = \mathbf{E}f(\theta) = 0$ and there is nothing to prove.

The proof of (i) is completed.

(ii) Using induction in $n$, we shall show that (i) yields (ii). For $n = 1$, the asertion (ii) is equivalent to (i). Assume that (ii) hold for $1, \ldots, n-1$. Let us prove (ii) for $n$. Notice that for given $X_1$, the sequence

$$Z_0 = 0, \qquad Z_1 = X_2, \ldots, Z_{n-1} = X_2 + \cdots + X_n$$

is a martingale sequence with differences satisfying

$$\mathbb{P}\{Z_k - Z_{k-1} \le b_{k+1}\} = 1, \qquad \mathbf{E}((Z_k - Z_{k-1})^2 | Z_1, \ldots, Z_{k-1}) \le \sigma_{k+1}^2$$

for $k = 1, \ldots, n-1$. Conditioning on $X_1$ and applying the induction assumption twice (for $n-1$ and 1), we have

$$\begin{aligned}
\mathbf{E}f(M_n) &= \mathbf{E}(f(X_1 + \cdots + X_n)|X_1) \\
&\le \mathbf{E}(f(X_1 + \theta_2 + \cdots + \theta_n)|X_1) \\
&= \mathbf{E}(f(X_1 + \theta_2 + \cdots + \theta_n)|\theta_2, \ldots, \theta_n) \\
&\le \mathbf{E}(f(\theta_1 + \theta_2 + \cdots + \theta_n)|\theta_2, \ldots, \theta_n) = \mathbf{E}f(T_n),
\end{aligned}$$

which completes the proof of (ii) and of the lemma. $\square$

LEMMA 4.5. *Let*

$$(4.9) \qquad x_1, \ldots, x_n \ge 0, \qquad a = (x_1 + \cdots + x_n)/n.$$

*Let $T_n = \theta_1 + \cdots + \theta_n$ be a sum of independent (eventually non-i.i.d.) Bernoulli random variables $\theta_k = \theta_k(x_k, 1)$. Let $S_n = \varepsilon_1 + \cdots + \varepsilon_n$ be a sum of $n$ independent copies of a Bernoulli random variable $\varepsilon = \varepsilon(a, 1)$. Let $f(x) = (x-t)_+^2$. Then, for any $t \in \mathbb{R}$, we have*

$$(4.10) \qquad\qquad\qquad \mathbf{E}f(T_n) \le \mathbf{E}f(S_n).$$

PROOF. Write

$$(4.11) \quad q_k = \mathbb{P}\{\theta_k = -x_k\} = 1/(1+x_k), \qquad p_k = \mathbb{P}\{\theta_k = 1\} = x_k/(1+x_k)$$

and notice that $\mathbb{P}\{\varepsilon = -a\} = 1/(1+a)$ and $\mathbb{P}\{\varepsilon = 1\} = a/(1+a)$.

We use well known properties of Schur convex functions [see Marshall and Olkin (1979)]. Recall that a vector $x = (x_1, \ldots, x_n) \in \mathbb{R}^n$ majorizes $y = (y_1, \ldots, y_n) \in \mathbb{R}^n$ (we use the notation $x \geq^* y$) if

$$x_{n:n} + \cdots + x_{k:n} \geq y_{n:n} + \cdots + y_{k:n} \qquad \text{for all } k = 1, \ldots, n,$$

where $x_{n:n} \geq \cdots \geq x_{1:n}$ is a decreasing rearrangement of the sequence $x_1, \ldots, x_n$. Notice that $x \geq^* y(x)$ for any $x \in \mathbb{R}^d$, where the vector $y(x) = (a, \ldots, a)$ has equal coordinates such that $a = (x_1 + \cdots + x_n)/n$.

A real valued function $g$ defined on an open subset $C \subset \mathbb{R}^d$ is called Schur concave if $x \geq^* y$ implies $g(x) \leq g(y)$. Assuming that $g$ has continuous partial derivatives such that

$$(4.12) \qquad \partial_j g - \partial_i g \geq 0, \qquad \text{when } x_i > x_j,$$

where $\partial_j = \partial/\partial x_j$, a result of Schur [see Schur (1923) and Ostrowski (1952)] says that $g$ is Schur concave in cases when the set $C$ is a symetric open convex set and $g$ is a symmetric function of its arguments. Notice that the result of Schur still holds if the set $C$ instead of the symmetry assumption satisfies: there exists a $z = (b, \ldots, b) \in \mathbb{R}^d$ such that the set $C - z$ is symmetric. Indeed, the majorization and (4.12) are preserved by a shift transformation of this kind.

Write $g(x) = g(x_1, \ldots, x_n) = \mathbf{E} f(T_n)$. Due to the result of Schur, to prove the inequality (4.10) it suffices to check that the function $g$ is a Schur concave function Notice, that as $C$ we can choose a sufficiently large open cube which, for a given $a$, contains the set

$$\{x \in \mathbb{R}^d : x_1 + \cdots + x_n = an, x_1, \ldots, x_n \geq 0\}.$$

Because the cube is open, we have to allow $x_k$ to assume negative values. We assume that $x_k \geq -1/3$. Now the probabilities defined by (4.11) can be negative, and in such cases we understand $\mathbf{E} w(\theta_k)$ as $\mathbf{E} w(\theta_k) = w(-x_k) q_k + w(1) p_k$.

Due to the symmetry of $g$ in its arguments, it suffices to check the condition (4.10) with $j = 1$ and $i = 2$. The inequality has to hold for all $t \in \mathbb{R}$. Therefore, conditioning on $\theta_3, \ldots, \theta_n$, it is easy to see that we can assume that

$$(4.13) \qquad g(x) = \mathbf{E} f(\theta_1 + \theta_2).$$

To simplify notation write $x_1 = \alpha$ and $x_2 = \beta$. Then

$$q_1 = \mathbb{P}\{\theta_1 = -\alpha\} = 1/(1 + \alpha), \qquad p_1 = \mathbb{P}\{\theta_1 = 1\} = \alpha/(1 + \alpha)$$

and

$$q_2 = \mathbb{P}\{\theta_2 = -\beta\} = 1/(1 + \beta), \qquad p_2 = \mathbb{P}\{\theta_2 = 1\} = \beta/(1 + \beta),$$

and we have to check that $\partial_\alpha g - \partial_\beta g \geq 0$ assuming that $\beta > \alpha$, where $\partial_\alpha = \partial/\partial\alpha$. For the function $g$ from (4.13) we have

$$(4.14) \qquad g = f(-\alpha - \beta)q_1 q_2 + f(1-\beta)p_1 q_2 + f(1-\alpha)q_1 p_2 + f(2)p_1 p_2.$$

We consider the following five cases separately:

   (i)  $t \leq -\alpha - \beta$;
  (ii)  $-\alpha - \beta \leq t \leq 1 - \beta$;
 (iii)  $1 - \beta \leq t \leq 1 - \alpha$;
 (iv)  $1 - \alpha \leq t \leq 2$;
  (v)  $t \geq 2$.

In the proof of (i)–(v) write $I(t) = \partial_\alpha g - \partial_\beta g$. The function $t \mapsto I(t)$ is a continuous function. We have to show that $I(t) \geq 0$.

*Case* (i). In this case $f(x) = (x - t)^2$ on the support of $\theta_1 + \theta_2$ and, therefore,

$$(4.15) \qquad g = \mathbf{E}(\theta_1 + \theta_2 - t)^2 = \alpha + \beta + t^2$$

and the inequality $I(t) \geq 0$ is just the equality $0 = 0$.

*Case* (ii). Now [cf. (4.14)]

$$g = (1 - \beta - t)^2 p_1 q_2 + (1 - \alpha - t)^2 q_1 p_2 + (2 - t)^2 p_1 p_2.$$

Adding and subtracting $(-\alpha - \beta - t)^2 q_1 q_2$ and using (4.15), we have

$$g = \alpha + \beta + t^2 - (\alpha + \beta + t)^2 q_1 q_2.$$

Using $\partial_\alpha q_1 = -q_1^2$ and $\partial_\beta q_1 = 0$, it is easy to find that

$$I(t) = (\alpha + \beta + t)^2 q_1^2 q_2^2 (\beta - \alpha) \geq 0,$$

which concludes the proof of case (ii).

*Case* (iii). We have [cf. (4.14)]

$$g = (t + \alpha - 1)^2 q_1 p_2 + (t - 2)^2 p_1 p_2.$$

Using $\partial_\alpha p_1 = q_1^2$ and $\partial_\beta p_1 = 0$, it is easy to see that

$$(4.16) \begin{aligned} I(t) &= 2(t + \alpha - 1)q_1 p_2 - (t + \alpha - 1)^2 q_1^2 p_2 + (t - 2)^2 q_1^2 p_2 \\ &\quad - (t + \alpha - 1)^2 q_1 q_2^2 - (t - 2)^2 p_1 q_2^2. \end{aligned}$$

The function $I(t) = At^2 + Bt + C$ is a quadratic function of $t$ with some $A$, $B$ and $C$. It is clear from (4.16) that $A = -q_1 q_2^2 - p_1 q_2^2 = -q_2^2$. This means that the function $t \mapsto I(t)\colon [1 - \beta, 1 - \alpha] \to \mathbb{R}$ is a concave function. Hence, $I(t) \geq 0$ will follow if we check the inequality at the endpoints of the interval. However, the inequality $I(1 - \beta) \geq 0$ is already established in (ii). The inequality $I(1 - \alpha) \geq 0$ is proved in case (iv).

*Case* (iv).   In this case $g = (t-2)^2 p_1 p_2$ and

$$I(t) = (t-2)^2(q_1^2 p_2 - p_1 q_2^2) = (t-2)^2 q_1 q_2 (\beta q_1 - \alpha q_2).$$

Hence, it suffices to check that $\beta q_1 - \alpha q_2 \geq 0$, which is equivalent to $(\beta - \alpha)(1 + \beta + \alpha) \geq 0$, which is obvious.

*Case* (v).   Now $g = 0$ and there is nothing to prove. The proof of the lemma is completed.   □

LEMMA 4.6.   *Let $f$ be one of the following functions:*

$$f(x) = (x-t)_+^2, \qquad t \in \mathbb{R};$$
$$f(x) = (x-t)_+^s, \qquad s > 2;$$
$$f(x) = \exp\{hx\}, \qquad h > 0.$$

(i) *Assume that a random variable $X$ satisfies*

$$(4.17) \qquad \mathbb{P}\{X \leq b\} = 1, \qquad \mathbf{E}X^2 \leq \sigma^2.$$

*Then we have $\mathbf{E}f(X) \leq \mathbf{E}f(\theta)$, where $\theta$ is a symmetric Bernoulli random variable $\theta = \theta(a^2, a)$ with $a = \max\{\sigma, b\}$.*

(ii) *Assume that a martingale $M_n$ satisfies (2.1), that is, $X_k \leq b_k$ and $s_k^2 \leq \sigma_k^2$ with probability 1, where $s_k^2$ are the conditional variances of the differences $X_k$. Let $T_n = \theta_1 + \cdots + \theta_n$ be a sum of independent symmetric Bernoulli random variables $\theta_k = \theta_k(a_k^2, a_k)$ with $a_k = \max\{\sigma_k, b_k\}$. Then $\mathbf{E}f(M_n) \leq \mathbf{E}f(T_n)$.*

PROOF.   Similar to the proof of Lemma 4.4, it suffices to establish (i).

Assume first that $\sigma \leq b$. By (i) of Lemma 4.4, we have $\mathbf{E}f(X) \leq \mathbf{E}f(\theta_0)$, where a Bernoulli random variable $\theta_0 = \theta_0(\sigma^2, b)$. The condition $\sigma \leq b$ implies

$$(4.18) \qquad \mathbb{P}\{\theta_0 \leq b\} = 1, \qquad \mathbf{E}\theta_0^2 \leq b^2.$$

In the view of (4.18) we can estimate the expectation $\mathbf{E}f(\theta_0)$ using (i) of Lemma 4.4. We get $\mathbf{E}f(\theta_0) \leq \mathbf{E}f(\theta)$ with a symmetric Bernoulli random variable $\theta = \theta(a^2, a)$ because $a = \max\{\sigma, b\} = b$, due to the assumption $\sigma \leq b$. Combining the inequalities, we obtain the desired $\mathbf{E}f(X) \leq \mathbf{E}f(\theta)$.

Assume now that $\sigma > b$. A random variables which satisfies (4.17), satisfies as well $\mathbb{P}\{X \leq \sigma\} = 1$ and $\mathbf{E}X^2 \leq \sigma^2$, and we can again apply (i) of Lemma 4.4, because now $a = \max\{\sigma, b\} = \sigma$.   □

LEMMA 4.7.   *Assume that for a function $Q(a; \sigma^2)$ the bound*

$$(4.19) \qquad \mathbb{P}\{M_n \geq an\} \leq Q^n(a; \sigma^2)$$

holds for all $n$ and all sums $M_n = \varepsilon_1 + \cdots + \varepsilon_n$ of i.i.d. Bernoulli random variables $\varepsilon_k = \varepsilon_k(\sigma^2, 1)$ so that the conditions of Hoeffding's Theorem 3 are fulfilled. Then we have

$$Q(a; \sigma^2) \geq H(aq + p; p), \qquad where \ p = \sigma^2/(1 + \sigma^2), \ q = 1 - p.$$

PROOF. We have $\mathbb{P}\{\varepsilon_k = -\sigma^2\} = q$, $\mathbb{P}\{\varepsilon_k = 1\} = p$ and

$$\mathbb{P}\{M_n \geq an\} = \mathbb{P}\{\theta_1 + \cdots + \theta_n \geq z\} \qquad \text{with } z = aq + p \text{ and } \theta_k = q\varepsilon_k + p.$$

The random variables $\theta_k$ are i.i.d. Bernoulli random variables such that $\mathbb{P}\{\theta_k = 0\} = q$ and $\mathbb{P}\{\theta_k = 1\} = p$. The inequality (4.19) implies

$$\log Q(a; \sigma^2) \geq \frac{1}{n} \log \mathbb{P}\{\theta_1 + \cdots + \theta_n \varepsilon \geq z\}.$$

Passing to the limit as $n \to \infty$ and using a well-known result on large deviations [see Bahadur (1971), Example 1.2], we get

$$\log Q(a; \sigma^2) \geq -f$$

with $f = z \log(z/p) + (1 - z) \log((1 - z)/(1 - p))$, for $p < z < 1$. Using the explicit formula (1.18) for $H$, it is clear that $\exp\{-f\} = H(a; \sigma^2)$, which proves $Q \geq H$ and the lemma. $\square$

LEMMA 4.8. *In Theorems 1.1–1.3, $\mathbb{P}\{S_n \geq x\}$ cannot replace $\mathbb{P}^\circ\{S_n \geq x\}$.*

PROOF. It suffices to prove the lemma in the case $n = 1$. Let $X$ be a random variable such that $\mathbb{P}\{X \leq 1\} = 1$, $\mathbf{E}X = 0$ and $\mathbf{E}X^2 \leq \sigma^2$. Let $\varepsilon = \varepsilon(\sigma^2, 1)$ be a Bernoulli random variable. To prove the lemma it suffices to check that

$$(4.20) \qquad\qquad \sup_{\mathcal{L}(X)} \mathbb{P}\{X \geq 0\}/\mathbb{P}\{\varepsilon \geq 0\} = \infty.$$

Taking $X = 0$ we have $\mathbb{P}\{X \geq 0\} = 1$. Using $\mathbb{P}\{\varepsilon \geq 0\} = \sigma^2/(1 + \sigma^2)$, we see that (4.20) is implied by the obvious $\sup_{\sigma^2 > 0}(1 + \sigma^2)/\sigma^2 = \infty$. $\square$

**Acknowledgments.** I thank I. Pinelis for pointing out the reference by Keilson and Gerber (1971). We thank N. Kalosha who performed the computer calculations that found the value of $c_1$.

VILNIUS INSTITUTE OF MATHEMATICS
AND INFORMATICS
AKADEMIJOS 4
232600 VILNIUS
LITHUANIA
E-MAIL: bentkus@ktl.mii.lt